\documentclass{article}

\usepackage{amsmath, amsfonts}
\usepackage{wrapfig}
\usepackage[dvips]{graphicx}
\usepackage[cp866]{inputenc}
\usepackage[english,russian]{babel}
\usepackage{amssymb, amsthm}
\usepackage{latexsym}
\newtheorem*{theorem*}{Threorem}
\newtheorem*{theorem*1}{Theorem \normalfont{on polyhedra}}

\newtheorem*{theorem*2}{Theorem on index}

\date{}

\title{\bf The Delone Peak \footnote{The author is very thankful to his friend,  professor Marjorie Senechal  for  spending  much  trouble on polishing an original text written  poor author's English. He also thanks A.B. Delone and V.S.Makarov for presenting photos.}}
\author{N.\,P.\,Dolbilin\\Steklov Mathematical Institute, \\ \small{dolbilin@mi.ras.ru} \footnote{ supported in part by RFBR project 01-08-00565-a
}}

\begin{document}

\maketitle

\renewcommand{\abstractname}{Abstract}
\renewcommand{\refname}{References}

\begin{abstract}
One of  beautifull mountains in the neighborhood of the Belukha,   the tallest  summit in the Altai  (in Siberia),  is called the Delone Peak in the  name of Boris Nikolaevich Delone (Delaunay), prominent mathematician. Here we are going to say about the life and science of this extremely suprising man \footnote{The paper is a seriously remade  and supplemented version of the author's paper "The Delone Peak" published in "Kvant", 1987, No (Russian)}.

\end{abstract}

\textbf{Childhood}

Boris Nikolaevich Delone (Delaunay) (March 15, 1890, Saint-Peterbourg - July 17, 1980, Moscow)
\footnote{In  childhood and in the first part of his life, publishing  papers in French,
  Boris Delone used the spelling \textit{Delaunay}. Beginning in the 1930s, he published his papers mostly in Russian with the Russian spelling of his name. When, later, his work was translated into  foreign languages his name was transliterated as
\textit{Delone}. Through out the  paper we will use the "Russian-English" spelling \textit{Delone} for Boris Nikolaevich  and  French spelling \textit{Delaunay} for  his ancestors. } was, above all, an outstanding mathematician, and mathematics  was certainly the main theme of his life. However, in a polyphonic musical piece it is hard to separate different themes. For people who knew Delone well, it is impossible to consider the  mathematical part of his life separately from his other also very strong, interests and activities.

Boris Delone's father,  Nikolay Delaunay\footnote{The Delaunay family originated in France: Boris' great grandfather Pierre Delaunay  was a medical doctor in Napoleon's army when it invaded Russia in 1812. He was sent from burnt-down Moscow to Simbirsk. After he release he returned to France. But, being in  love with a young  beautifull gentlewoman from a noble Russian family, Pierre   returned to Russia, married  and  remained  there for the rest of his life. Reports that the Commandant of the Bastille in Paris,  Marquis De Launay, an early victim of the French Revolution  killed by rebels on July 14, 1789, was Boris Nikolaevich's ancester turned out to be wrong. According to later investigations  by professor of astronomy Alexander Sharov (B. Delone's son-in-law), the Marquis  belonged to another family.},  was a well-known  professor of mechanics, an author of university textbooks. For Boris' biography, it turned out to be important that in the early 1900's, while his father was a professor of Warsaw university\footnote{At that time the territory of Poland was a part of the Russian empire.} he was a  close friend of Georgy Voronoi. Later Boris Delone used to recall that Voronoi  visited his home and the young teenager Boris often attended evening  conversations between his dad and Voronoi. 
 
 Boris Delone had an excellent, many-facet  education. His studies in music were quite solid: he played many pieces  of Bach and Mozart, all the sonatas of  Beethoven, and composed many himself. A music teacher recommended that the musically gifted boy enter a conservatory to study  composition. However, his drawing teacher insisted that he continue his education in an art academy!

While parents and teachers tried to determine ``the  future'' of Boris, the boy painted landscapes and played soccer, reproduced (in pencil) "The Last Supper" by Leonardo, and climbed trees with a younger sister on his shoulders (``for physical load'', Delone explained later). At the same time Boris turned his room into a real physics laboratory in which he made many of the devices himself. Later he often recalled proudly  small technical tricks  which allowed him, for example, to obtain, in a Leyden jar he made, a ``so-o-o huge spark''. Being keen on astronomy, Boris made  a telescope for which  he also polished a mirror. Later, each time he mentioned this episode, he said that "polishing a mirror  from bronze was a stupid idea: it was very labor-consuming and grew dim fast".

 It is impossible not to mention one more event in Boris Delone's youth. His father
was a friend of a famous Russian scientist Nikolay Egorovich Zhukovski, the "grandfather of Russian aviation''.  Under his influence in 1907 Nikolay Delaunay organized the first gliding circle in Russia in Kiev,  after the family moved from Warsaw. The 17 year old Boris took an active part in this circle. During the next two years he constructed 5 gliders and flew on them. Needless to say we are talking here  about flights  of several tens of meters. Nevertheless, this circle played a notable role in development of aviation in Russia and in the world. In particular, among the participants was Boris's friend Igor Sikorski, also 17.
Later, in 1913, Igor Sikorski constructed a huge (for that time) four-engine aircraft. After the Bolshevik  revolution in 1917, Sikorski emigrated to the United States, where very soon he organized a famous helicopter company. One of his helicopters was the first to cross the Atlantic.

Before discussing the place of mathematics in Boris's childhood, one should mention  one more passion that appeared at age of 12-15  and continued throughout his life: a passion for mountaineering and hiking. This passion was began and grew in early 1900's in the Swiss Alps, where Boris' family used to go for summer vacations. Later Delone talked about his climbing and hiking around Mont Ros and the Matterhorn above Zermat.

Boris Delone's mathematical talent was evident very early. By age 12 Boris was familiar with the elements of Calculus. At 13  he started his own investigations in Algebra and Number Theory. His father encouraged Boris in his mathematics studies and in 1904 took him to the  International Congress of Mathematicians in Heidelberg. There Delone attended talks by Hilbert and Minkowski. Voronoi too took an active part in the congress, and got to know Minkowski. His conversation with Minkowski likely inspired Voronoi's tempest activity in the geometry of numbers in the last years of his life.

 \medskip
 \textbf{Univerity years}

 \smallskip
  Voronoi's work strongly influenced Delone's, though Boris had no direct scientific contacts with him. Voronoi died suddenly at the age of 40 in 1908, exactly the same year that Delone entered the physics and mathematics department of Kiev university.

Several extremely talented young people entered the department at about that time,  among them  Otto Yu.Schmidt (1891-1956)  and Nikolay Chebotarev (1894-1947). Chebotarev became an outstanding algebraist. Schmidt also became a prominent algebraist. He wrote a well-known book ``Abstract theory of groups'', founded a famous Moscow algebraic school, and became a full member of the USSR Academy of Sciences.
 Moreover, the name of Schmidt became legendary in the Soviet Union as  an explorer of the Arctic Ocean and a leader of a scientific expedition on the ice-breaker "Chelyuskin" in 1934.

 Boris Delone, Schmidt, Chebotarev, and several others formed the core of the seminar led by Dmitrii Grave. This seminar on algebra and the algebraic theory of numbers, legendary in the  history of Russian mathematics, determined Delone's field of studies for many years. Delone's first student work,  ``Link between the theory of ideals and the Galois theory'', was awarded  the ``Great golden medal'' of  Kiev university. Delone's  first  published paper, ``On determining algebraic domain by means of congruence'', was dedicated to a new proof of a famous theorem of Kronecker on absolutely abelian fields.

\medskip
\textbf{Cubic diophantine equations}

\smallskip

 At about this time Delone began his investigation of the theory of Diophantine equations of degree 3 with two variables. This turned out to be the summit of  his entire mathematical career.

 Diophantine equations are
 of the  form
 $$
 p(x_1,\ldots, x_n )=0, \eqno(1)
 $$
 where $p(x_1,\ldots , x_n)=0$ is a polynomial with integer coefficients, for which one needs to find integer (or sometimes rational) solutions. For example, the equation $x^2+y^2=z^2$ has infinitely many integer solutions, the well-known Pythagorean triples (3,4,5), (5,12,13), and so on. These triples correspond to right-angled triangles with integer sides.

 In Hilbert's famous list 23 problems the 10th problem asked: is there an algorithm to determine, from the coefficients of equation (1), whether or not solutions exist? Today, after the work by Yu. Matiyasevich, it is well-known that there is no such algorithm. In fact, the negative answer was expected by all experts on diophantine equations. They were very well aware of how hard each new result in this area was won.

 The simplest form of diophantine equation is a linear equation with two variables:
  $$
  ax+by=1,
  $$
  where $a$ and $b$ are integers. This equation was investigated by the Indian mathematician Ariabkhata in 5th - 6th centuries.
 Diophantine equations of degree 2 turned out to be more difficult. They were studied  by Euler, Lagrange (a full theory of Pell's equation $ax^2+y^2=1$, where $a$ is a non-square number), Gauss and  other outstanding mathematicians of 19th century.
 As for diophantine equations of degree 3 or greater, a very deep result was obtained by the Norwegian mathematician A. Thue in 1908. Thue proved that diophantine equations of the following type
 $$ f(x,y)=c,
 $$
 where $f(x,y)$ is a homogenous polynomial of degree 3 or greater, has at most a finite number of (integer) solutions,  if the polynomial $f(x,y)$ cannot be decomposed into a product of polynomials of  degree 2.
This fact contrasts with  quadratic Pell's equations which have an infinite number of solutions.
One should note that Thue's theorem does not give any tools for finding these solutions.

It may be fortunate that the student Delone did not know Thue's work. Delone decided to investigate a longstanding problem of cubic diophantine equations equations of the form
 $$ f(x,y)=1,
 $$
 where $f(x,y)$ is a homogenous form of degree 3 with negative discriminant. (Delone understood that if any Diophantine equations of degree could have infinitely many solutions, it would be this class.) First, he studied a cubic analogue of  Pell's equation:
$$
x^3q +y^3=1, \eqno(2)
$$
where $q$ is an integer but not a cube. (Otherwise, a form $x^3q +y^3$ would be decomposable and the problem would reduce to solving diophantine equations of degree 1 and 2.)

Equation (1) has always a trivial solution (0; 1). The problem is to find other solutions, if any.
 Delone introduced
into consideration a ring $\sigma$ of algebraic numbers $z\sqrt[3]{q^2} + x\sqrt[3]{q} + y$, where $(z,x,y)$ are integer triples. It turns out that if $(x,y)$ is a solution (in integers) of equation (2) then
 a number $\varepsilon= x\sqrt[3]{q} + y$ belongs to the ring $\sigma$ together with its reciprocal $\varepsilon ^{-1}$ . Such elements of the ring $\sigma$ are called  \textit{units}. Thus, any solution of (2) corresponds to a unit of $\sigma$. Moreover, this unit $\varepsilon$ is a bi-term unit since one term,  $z\sqrt[3]{q^2}$, is equal to $0$.
It follows from a theorem of Dirichelt that the bi-term unit $\varepsilon\in \sigma$ is equal to a power $m$ of some so-called \textit{basic unit} $\varepsilon_0$ of $\sigma$: $\varepsilon=\varepsilon_0^m$. Indeed, it was also known that the degree $m>0$ if $0< \varepsilon<1 $.

Further, Delone proved that if the basic unit $\varepsilon_0$ is bi-term itself then no positive power $m$ of $\varepsilon_0$, except $m=1$, is equal to a bi-term unit. After that he investigated the case $\varepsilon=\varepsilon_0^m$, where $\varepsilon_0$ ia a three-term basic unit. By means of sparkling wit Delone proved that a bi-term unit of sort $\varepsilon_0^m$ can be just a bi-term basic unit $\varepsilon_0$ itself.

Thus, Delone obtained a final result: the cubic analogue of Pell's equation (2), besides a trivial solution (0; 1), has at most one more non-trivial solution. In order to get this solution one needs to find a basic unit. Here one needs to emphasize that about 20 years earlier, Voronoi had constructed an algorithm for computing basic units in the algebraic ring $\sigma$. If this unit is bi-term of the sort $x\sqrt[3]{q} + y$ then $(x,y)$ is the only non trivial solution of equation (2). If this is not the case then there are no non-trivial solutions at all.

After this strong success Delone proceeded to a more general problem:
$$f(x,y)=1, \eqno(3)$$
where $f(x,y)=x^3 + ax^2y +bxy^2+cy^3$ is a cubic form in two variables with a negative discriminant.

% I have deleted the definition of negative discriminant because you used
% the term already (without defining it), your definition is not fully general,
% and the reader can look it up on line.

In this more general case Delone again reduced the problem to an investigation of bi-term units of sort
$x+y\rho$ in a ring $\sigma$ of algebraic numbers $x+y\rho +z\rho^2$. By means of his so-called "algorithm of ascent" Delone proved the following fundamental theorem:

\textit{In general case equation (3) has at most 3 integer solutions; in two concrete cases it has exactly 4 solutions; one concrete equation has sharply 5 solutions. No equations of sort (2) have more than 5 integer solutions.}

% what do you mean by "sharply"? Is this different from "exactly"?

Moreover, by means of the algorithm of ascent, Delone  was able to find all solutions for each concrete equation of form (3).
For example, the equation $x^3 - xy^2+y^3=1$ has the 5 solutions (1; 1), (1; 0), (0; 1), (-1; 1), (4; -3).
The fact that the equation has no other solutions follows from the theorem by Delone mentioned above.

 However, Delone failed to rigorously prove that this method of obtaining solutions is workable for any equation. The matter was that then there was no upper bound expressed in terms of coefficients of the equation (3). Therefore, it was impossible to indicate the point when the   algorithm of ascent had already given all solutions or not in a general case. The existence of such a upper bound was obtained in 1960's by A. Baker. This work of Baker has been awarded of the Fields medal.

Finally, one should emphasize that after classical results by Euler, Lagrange, Gauss et al on quadratic diophantine equations, works by Delone represented a serious breakthrough in the theory of cubic equations and remained unsurpassed until the late 1960's (A.Baker). Boris Delone himself evaluated his investigations on cubic equations with two unknowns as the best  in all his scientific work.

\medskip

\textbf{Life in Kiev, move to Saint-Peterbourg (Leningrad)}

\smallskip
The success was a result, as Delone used to say, of thousands of hours of very intensive work. One should say that Delone managed to win a very complicate problem regardless extremely hard political and living conditions which existed in Kiev. During WW I Kiev was captured by German troops. Therefore, in 1915-1916 the Delone family (together with other university professors) was forced to leave Kiev for Saratov. Nikolai Chebotarev went to Saratov too. Over there Boris Delone strongly influenced  Nikolai, who was 4 years younger.

When the Delones returned to Kiev, in Russian territory, World War I was giving way to a more cruel and bloody civil war (1918-1922). Kiev was at the core of this war too. One political and military power was superceeded by another permanently. German troops left the city to the Red army (bol'shevics); in their turn the reds left Kiev to the White army (protectors of the Russian tsarist monarchy). The whites were pushed out again by the reds. The reds were replaced by Petlura's squadrons (a sort of ukranian national movement), which  were replaced by Polish troops, then by the greens, the greens by the yellow-blues and so on.

Later B.N. Delone laughed about one dramatic and illustrative story from his family's life. B.N. had a younger brother Alexander (nicknamed Alik) who was an officer in Denikin's army\footnote{General Denikin was the White army comander, one of the most prominent activists of the White movement in Russia.}.
Once, during a sudden capture of Kiev by red troops,  Alik was forced to escape from Kiev. For this he exchanged his military uniform for civil clothes. The military form was left in
a two door wardrobe in his father's apartment, where Boris also lived. The next day a red marine patrol came to them for inspection. During a search of the apartment a commander of the patrol approached the wardrobe  and abruptly opened the right door. Suddenly a big dol fell from an upper shelf onto the floor and started crying loudly. This doll had been brought from Paris for someone in the family a few years before. At that time a crying doll was a real wonder. It surprised all the patrols so much that the commandant forgot to open the other door of the wardrobe behind which the white officer uniform hung... The time was very cruel and ruthless, shooting was an ordinary matter.

In this period Delone worked successively as a math teacher in a gymnasium, an educator at Kiev university, and an associate professor at Kiev polytechnical university. In 1920 he presented his work on cubic diophantine equations as a doctoral dissertation to Saint-Peterbourg (at that time named Petrograd) university. Andrei Markov (of the Markov processes) was a leader of a commission which considered the dissertation and highly evaluated the thesis. In 1922 Delone was invited to be a professor of Saint-Peterbourg university. The university hosted a famous Saint-Peterbourg  school of number theory which had been begun by Euler and flourished in the late XIX century under P. Chebyshev's. A. Markov, G. Voronoi, A. Korkin, E. Zolotarev, A. Lyapunov and others formed the core of this school.

 In due time, according to the traditions of this famous number theoretical school, the outstanding member G. Voronoi clothed a very geometrical work  on parallelohedra in analytical clothes. In contrast, Boris Delone, possessing sure geometrical gift, geometrized numerous algebraic works, including the Voronoi algorithm for computing the basic unit in a ring which belongs to a field of degree 3 with a negative discriminant. After that he constructed a very nice geometric theory of cubic binary forms. In these and other works Delone interpreted the ring of cubic irrationalities as an integer lattice with a natural multiplication rule.

\medskip

 \textbf{The empty sphere method}

 \smallskip

Delone's study of algebraic problems with geometric tools  continued into the late 1950s. Here we are not going to focus any longer on this part of his activity. It is interesting that already in 1920s he wrote a paper on the  empty sphere. A tendency to clarity and transparency  was always peculiar to Delone. He gave an extremely  clear description of the geometry of a tiling dual to the well-known Dirichlet tiling. Delone  called this dual tiling a $L$-\textit{tiling}, one can only guess why. Likely, it was  in honor of the first part of Dirichlet's surname (his full name was Johan Peter Gustav \textit{Lejeuene Dirichlet}) \footnote{Dirichlet, (1805-1859, used such 2- and 3-dimensional tessellations in his studies on quadratic forms.}

 For this,
 Delone introduced an \textit{(r,R)-set} $X$, where $r,R$ are positive numbers. Today such sets are also called  \textit{Delone sets}, or  \textit{separated nets}. By definition, such a set fulfils:
\newline (1) any open ball with radius $r$ contains at most one point of $x$;
\newline (2) any closed ball with radius $R$ contains at least one point of $x$.

Next,  Delone displaced into space a small ball free of points of $X$. Such an empty ball exists due to the property (1). Then, leaving its center fixed,  he expands the ball until it touches at least one point $x$ of $X$. This event must happen due to point (2). Moreover, it has to occur while the radius still remains smaller $R$.
 Next, keeping $x$ on its boundary, Delone expands  ball (still empty inside) away from $x$. Sooner or later a new point $x'$ of $X$ appears on the boundary of the growing ball. At the next stage Delone increases the ball further, keeping the points $x$ and $y$ on its boundary. This continues until the boundary contains a full set of independent points. If it is to remain empty, the ball cannot be increased further. The radius of this locally maximal ball $B$ cannot exceed $R$ (due to condition (2)).
 By taking the convex hull of the points of $X$ on the boundary of $B$ Delone gets a convex inscribed polyhedron he called an $L$-\textit{solid}. The last crucial point is that the set of all $L$-solids forms a face-to-face tiling of space. Delone called it an  $L$-\textit{tiling}.

 Needless to say, in the 1960-70s Delone and we, his students, working with these tilings, did not suspect that very soon they would be called \textit{Delaunay triangulations}. I remember how proud Delone was  of the crystallographic terminology \textit{Delone symbol} and \textit{Delone sort}. However, the use of these terms, named in honor of their inventor, cannot compare in use and dissemination with Delone triangulations.

Unfortunately this term came to Russia from the West after Delone's death. Here one should note a contribution Delone and Coxeter made together for this term to enter the mathematical language.

It was rather curious story. One day, sometime in the late 1950s, Delone saw a paper by H. S. M.  Coxeter in which Coxeter had introduced and used a Delone tiling with no reference to Delone's work. Delone wrote Coxeter a letter in which, as Delone told it, he informed Coxeter that he had studied these tilings in the early 1920s. Delone asked Professor Steklov \footnote{Vladimir Andreevich Steklov (1864 - 1926) was a prominent Russian mathematician, academician  and, in the 1920's, organizer and the first director of the Institute of Mathematics and Physics of Russian Academy of sciences. Later, in 1934, this institute was divided into two institutes; one of them is Steklov Mathematical institute of Russian academy of sciences.}, who went in 1924 to an international congress of mathematicians held in Toronto, to deliver
 Delone's work. This paper was published in the proceedings of the congress\footnote{Sur la sph\'ere vide. Proceedings of the International mathematical congress held in Toronto, August 11 - 16, 1924, V. 1, Toronto, Univ. of Toronto press, 1928, p. 695 - 700.}. In addition, in his letter Delone mentioned that the congress was held in the "Coxeter's native city," Toronto. Coxeter's reply was, Delone said, very polite. He found, in the materials of Toronto congress, some mention on Delone's work. Coxeter asked Delone to forgive him because in 1924 he, Coxeter, was so young that he ``walked under table in short pants''. At the same time Coxeter wrote a letter about this to Prof. C.A.Rogers in Cambridge who was completing his book, ``Packings and Coverings''. Rogers inserted a mention of the Delaunay tiling. Very likely the term of Delaunay tilings (in Computational Geometry - triangulations) first appeared in C. A. Rogers' book.

 In 1993, I had a chance to meet Professor Coxeter in person and I asked him about this story. On the whole Professor Coxeter  confirmed it, except for ``walking under table in short pants''. In 1924 Coxeter (1907 - 2003) was already 17 years old. The phrase ``walking under the table'' was typical of Delone's expressive and humorous style.

\medskip
\textbf{Parallelohedra and the Geometry of Positive Quadratic Forms}

\smallskip
The work on the empty sphere began a new line in scientific work for Delone. If earlier his geometric gift showed in the geometrization of algebraic and number theory, in the late 1920s geometry itself became the main object of his interest. In 1929 Delone completed a large paper in which he derived all combinatorial types of 4-dimensional parallelohedra\footnote{Indeed, Delone found 51 parallelohedra and missed one which later, in 1969, was discovered by his student M.I.Stogrin.} In this paper Delone continued the classical work of Fedorov (all 5 types of 3-dimensional parallelohedra) and Voronoi (general theory of Voronoi parallelohedra for arbitrary dimension, all 3 primitive 4-dimensional parallelohedra). A parallelohedron of dimension $d$ is a convex $d$-polyhedron which tiles euclidean $d$-space, i.e. paves space without gaps, in face-to-face way. This concept was introduced by the crystallographer Evgraph Stepanovich Fedorov.

A very important special case of parallelohedra is the Voronoi domains for points of integer lattices in space. These parallelohedra  are called now Voronoi parallelohedra as Voronoi, in  his two last outstanding memoirs of 1908, invented and developed a general theory of such parallelohedra. By the way, these memoirs contain certain very deep ideas which still remain unknown to modern mathematicians. In particular, a well-known idea on the lift of a Voronoi tiling to a paraboloid of revolution, discovered by computational geometers about 30 years ago, was introduced and actively used in Voronoi's memoirs. Voronoi had a strategic plan to investigate first of all Voronoi parallelohedra and constructed an exhaustive theory. In parallel, he posed a conjecture: any parallelohedron is affinely equivalent to some Voronoi parallelohedron. This conjecture was proved by Voronoi for the most generic class of parallelohedra - primitive parallelohedra. (Recall that a parallelohedron is \textit{primitive} if at any vertex of a tiling by such parallelohedron the minimum  number of cells (i.e. $d+1$) meet. Later O.Zhitomirski proved affine equivalence of some Voronoi parallelohedron of any parallelohedron which is just \textit{primitive  at $(d-2)$-faces}, i.e. the minimum number (three) of cells meet at any   $(d-2)$-face   \footnote{O.K. Zhitomirski was a student of Delone}.
%% THE SENTENCE IS INCOMPLETE

Since this result the Voronoi conjecture has remained unsolved for over 70 years.

 In his work on 4-parallelohedra, Delone proved that \textit{any 4-dimensional parallelohedron is affinely equivalent to some Voronoi parallelohedron}. For dimension 5 and higher the answer to the Voronoi conjecture remains unknown.

At this time Delone started the further development of a theory which appeared in works by Gauss and was further developed by Voronoi. Delone introduced a term for this field: \textit{geometry of positive quadratic forms}. The main idea of the geometry of positive quadratic forms (PQF geometry) is as follows.

Given an integer lattice $\Lambda \subset \mathbb{E}^d$ of rank $d$, i.e. the set of all points $ \textbf{x}=(x_1,\ldots , x_d)$  all of whose co-ordinates $x_i$ are integral with respect to some vector basis ${\cal E}=\{\textbf{e}_1,\ldots , \textbf{e}_d\}$. A $(d\times d)$-matrix $A=(a_{ij})$, where $a_{ij}=(\textbf{e}_i, \textbf{e}_j)$, is symmetric matrix called a Gram matrix. A quadratic form $\textbf{x}A \textbf{x}^{\small T}= \sum_{ij}a_{ij}x_ix_j$ is a positively definite form. Each quadratic form $f=\sum_{ij}a_{ij}x_ix_j$ with $d$ variables is assigned to a  point $f\in \mathbb{E}^N$, where $N=\frac{d(d+1)}{2}$. Positive forms occupy in $\mathbb{E}^N$ an open convex cone $K$.
Since a lattice $\Lambda$ has infinitely many bases ${\cal E}_1, {\cal E}_2, \ldots $, each lattice $\Lambda$ is associated with infinitely many equivalent points $f_1, f_2,\ldots \in {\cal E}$. The cone $K$ is divided into classes of equivalent forms, i.e. forms corresponding to the same lattice $\Lambda$. Therefore, a primary task of PQF geometry is the reduction theory which chooses from each class an appropriate representative, the so-called \textit{reduced form}.

The merit of Delone's PQF geometry is  that he looked at the basic problems of the PQF theory from a unified geometric viewpoint. The empty sphere method lies at the core of  PQF geometry, as does the interpretation of a $d$-dimensional lattice as a point in a  convex cone  in a high dimensional space $\mathbb{E}^{\frac{d(d+1}{2}}$.
On this geometric foundation, Delone was  able to systematize and enrich the results of his predecessors.

We mention here just one result relevant to the optimal covering of space by balls.
Given a lattice $\Lambda\subset \mathbb{E}^d$, take the set of all equal balls
$$\{B_{\textbf{x}}(\rho)\, |\, \textbf{x}\in \Lambda\} \eqno(4)$$
centered at points $\textbf{x}$ of $\Lambda$ with radius $\rho$. One supposes here that $\rho$ is a minimal covering radius, i.e. the set (4) of balls covers $d$-space, but balls of smaller radius do not.  The main goal of the theory is to find a lattice (with a unit volume of the fundamental parallelepiped) giving the thinnest covering of space by equal balls. Due to efforts of Voronoi, Delone and others, the problem of finding of the thinnest lattice covering is now reduced to the following. Voronoi  found that the reduction domain in the cone $K$ with apex $O$ is divided into finite number of so-called \textit{$V$-type (Voronoi) domains} $V_1, V_2\ldots $. Each $V$-type domain $V_i$ is a polyhedral cone with the apex $O$. Each inner point $f\in V_i$ in the domain $V_i$ represents a lattice $\Lambda$ such that a Voronoi parallelohedron in the tiling $V(\Lambda )$ is primitive and the type of Voronoi tiling $V(\Lambda )$ remains unchangeable whi
 le $\Lambda$ stays inside of $V_i$.

The Delone tiling $D(\Lambda )$ for this lattice is dual to the $V(\Lambda )$. Therefore, the $L$-type (or Delone type) domain coincides with the $V$-type domain. From the definition of a Delone cell it immediately follows that the covering radius for a lattice $\Lambda$ is equal to the biggest radius of an empty sphere circumscribed around simplices
 in the Delone tiling.  Take a simplex $S({\cal E})$, whose all vertices have integer coordinates with respect to a certain basis $\cal E$, and  let its circumradius be $R(S, {\cal E})$. It turns out that a $R_0$-level surface $F(R_0)\subset \mathbb{E}^{\frac{d(d+1)}{2}}$ consisting of all bases $\cal E$ such that all circumradii $R(S, {\cal E})=R_0$ is \textit{convex}\footnote{B.N.~Delone, N.P.~Dolbilin, S.S.~Ryshkov and M.I.~Stogrin. A new construction in the theory of lattice coverings of an $n$-dimensional space by equal spheres. {\em Math.\ USSR Izv.} {\bf 4} (1970), 293--302.}. Due to this important fact, in each $L$-domain there exists at most one locally thinnest lattice.  Therefore, the problem of the optimal lattice is reduced to a concrete optimization problem for any given $L$-type domain. After that by comparing the maxima over all finitely many $L$-type domains one can get the thinnest lattice covering of space by equal balls.

Among  Delone's other results  geometrizing the algebraic theory of numbers we mention a brilliant geometric description of the Voronoi algorithm for computing a basic unit of a ring contained in the field of the third degree with a negative discriminant\footnote{B.N.Delone. Saint-Peterbourg school of the number theory, 2005.}. This was not by chance. As we remember, a  basic unit of an algebraic integer ring and its computation played an important role in Delone's research on diophantine cubic equations.

\medskip
\textbf{Moscow period}

\smallskip
From 1922 to 1934 Delone stayed in Leningrad (now Saint-Peterbourg)and taught in Leningrad university. Due to his outstanding pedagogical and artistic talents Delone became famous as an outstanding, remarkable professor. In 1930-1934 he held the algebra and number theory chair of Leningrad university. In 1929 Delone was elected to the USSR Academy of Sciences as a corresponding member. He joined Institute of Mathematics and Physics just founded by Steklov.
In 1934 the Institute was divided into the Steklov Mathematical Institute and the Lebedev Physics Institute. Since that time Delone's career has been tied  to
the Steklov institute of the Academy of Sciences. That same year, 1934, the USSR Academy of Sciences' headquarter moved to Moscow. The Steklov Institute was subdivided into two parts and one of them moved to Moscow too.

Delone joined Moscow branch of the Steklov Institute where he was  the  head of the algebra department (1945 - 1960) and later (1960 - 1980)  the head of the geometry department. Besides his research at the institute,  Delone was doing intensive pedagogical work, which continued into the early 1960s,
 teaching at Moscow university (professor, 1935-1958, and the head of the higher geometry chair, 1935 - 1943).  His  long teaching experience resulted in a two-volume textbook on analytical geometry\footnote{B.N. Delone, D.E. Raikov, Analytical Geometry, V. I, (1948), pp. 456, V. II, (1949), pp. 516 (in Russian).}. Indeed, this enormous book can be  thought of as a geometric encyclopaedia. The famous topologist Paul Alexandrov said that due to the beauty of its geometric ideas, Delone's book has no equal.
Delone was also an organizer and the first head of the mathematics chair of the newly inaugurated Moscow Institute of Physics and Technology (1947 - 1961).

During his Moscow period, Delone did research mainly on PQF geometry, the geometrization of diophantine cubic equations, and application of mathematics  to crystallography. In 1937-38 he published a large memoir (in two parts) on PQF\footnote{B.N.Delone. Geometry of positive quadratic forms, P.I, Russian surveys, \textbf{3}, (1937), p.16-62; P.II, Russian surveys, \textbf{4}, (1938), p.102-164, (in Russian)}. This paper gave a name to the whole field of geometry of numbers, as we had already said, the name of \textit{PQF geometry}. Geometric methods of studying cubic diophantine equations and cubic irrationalities were summarized in a monograph written jointly with his closest student, the outstanding algebraist D.K. Faddeev\footnote{Delone B.N., Faddeev D.K.. The theory of irrationalities of the third degree, M. - L., (1940), pp. 340 (in Russian), English translation: Translations of Mathematical Monographs V. 10, Providence, Rhode Island, AMS, (1964), pp. 513.}.

Among his crystallographic works one must mention Delone's discovery of the so-called \textit{Delone sorts of lattices}.  Delone defines two lattices $\Lambda$ and $\Lambda '$ to be of the same sort if their Voronoi tilings $V(\Lambda)$ and $V(\Lambda ')$ are combinatorially isomorphic, the symmetry groups $Sym (V(\Lambda ))$ and $Sym (V(\Lambda '))$ are isomorphic, and the isomorphisms commute. Delone classified all 3-lattices into 24 Delone sorts. This classification involves not only the coincidence of lattice group properties,  but also
the geometry of the regular tilings corresponding to a given lattice. Therefore Delone's classification  was a subsequent, more detailed development of Bravais classification of lattices, and proved itself a very natural and helpful for crystallography,

\medskip
\textbf{Isohedral Tilings and Local Theory}

\smallskip
In the late 1950s B.Delone showed a strong interest in the theory of isohedral tilings. A tiling is called \textit{isohedral} \footnote{Delone himself called such tilings \textit{regular}.}
if its symmetry
 group operates on the set of tiles in transitive way, i.e. for any pair of tiles $P$ and $P'$ of a tiling there a symmetry of the tiling which maps $P$ to $P'$.
An isohedral tiling is assumed to be face-to-face and
the tiles convex. The cells of such a tiling are called \textit{stereohedra}\footnote{The term a stereohedron is due to E. Fedorov. A parallelohedron is a special case of stereohedron.}. Due to the Schoenflies-Bieberbach theorem\footnote{This principal theorem on the structure of crystallographic groups was an answer to  Hilbert's 18th problem. In fact this problem consisted of two problems. One  was to prove that for any dimension $d$  any crystallographic group operating in Euclidean $d$-space has a translational subgroup of finite index.}, any isohedral tiling with stereohedra is partitioned into several "lattices" of stereohedra. The set of stereohedra in one lattice is invariant and transitive with respect to the translation subgroup.

In contrast to the theory of parallelohedra, well-developed by H.Minkowski, G.Voronoi, B.Delone, et al, the theory of stereohedra remained undeveloped. Jointly with N.Sandakova Delone found the upper bound for the number of facets of a $d$-dimensional stereohedron in the given tiling:
$$ f_{d-1}\leq 2(2^d-1)+(h-1)2^d, \eqno(5)
$$
where $h$ denotes the number of lattices of stereohedra in the tiling. Since $h$ is upper bounded for the given dimension $d$ (for example, for all $d\geq 10$ $h\leq 2^dd!$) (5) implies a universal upper bound for the number of facets in any $d$-dimensional stereohedra. Thus, Delone generalized a celebrated upper bound by Minkowski for the number of facets in a parallelohedron,
$$
f_{d-1}\leq 2(2^d-1). \eqno(6)
$$
Though the upper bound (6) is obtained from (5) by putting $h=1$, it is true that, in contrast to (5) which is non-refinable, the upper bound (5) is very rough for $h>1$.
% Which do you mean, 5 or 6?
Nevertheless, from this upper bound the upper bound for the number of combinatorial types of isohedral tilings of euclidean space follows for any given dimension.
This was a key idea in constructing a general theory of euclidean stereohedra.

One should mention also the so-called \textit{local theory} of isohedral tilings and regular point sets. An isohedral tiling (or regular point set = crystallographic orbit) is an appropriate mathematical model of a crystalline structure. The problem is that these mathematical concepts are defined by the notion of symmetry group yet a real crystalline structure is formed without any idea of a group: it is just a result of the interaction of nearby atoms. Assume we are dealing with a solution consisting of atoms of one sort. Under appropriate physical conditions (concentration, temperature, pressure) the crystallization process starts. Each atom  tries to surround itself with other atoms in exactly the same way. So, it is likely that the source of global order in a crystalline structure is caused by the local identity of the structure in the neighborhood of each atom.
Something like that can be watched in an army. When soldiers are stepping in rank each soldier aligns himself not with all rest but just with the fourth one on the right. Thus a \textit{local rule}  provides a \textit{global order} in a rank.

Early on, Delone and his students suggested searching for mathematical theorems explaining the link between local identity  and global order in crystals\footnote{B.N.Delone, N.P.Dolbilin, M.I.Stogrin, R.V.Galiulin. A local criterion for regularity
of a system of points. Soviet Math. Dokl. 17 (1976), 319-322}. Later a more complete local theory, whose main goal was searching for and describing local conditions in discrete structures that provide  global crystallographic order, was developed in a series of works by M.I.Stogrin, N.P. Dolbilin, and others.

With his research, Delone actively attracted young people and influenced their progress. Among his students there are several outstanding scholars who themselves founded their own schools. They include, most notably, the geometer academician Alexander Alexandrov, and the algebraists academician Igor Shafarevich and corresponding member of USSR Academy of sciences Dmitrii Faddeev. Among others one should mention the names of Tartakovskii V.A., Zhitomirskii O.K., Ryshkov S.S..

\medskip
\textbf{Not only science}.

\smallskip
The Delone's scientific work is varied, valuable, and held in very high respect. However, the extreme brilliance of Delone's personality consists in that, besides outstanding scientific talent, he possessed  other  talents also and, what is most amazing, he was able to demonstrate explicitly  these talents in his life.
Likely, the ability to live as independently as possible, according to his own principles and goals, among which service to mathematics stood in first place, was Delone's  main talent.
I remember well how one gray day in March 1970 the great Andrei Kolmogorov came to Delone's small study in the Steklov to congratulate him on the occasion of his 80th birthday. Kolmogorov presented a gift (it was a rolled up gravure showing, in the  old Dutch style, a marine landscape with war sailing vessels). During their talk Kolmogorov mentioned  how much he envied Delone  that he (Delone) had been able to avoid the temptations of seductive offers of this or that high position in scientific administration or something like that, and to organize his life by his own principles and  enjoy daily life. It was really true to fact...

Delone was a fantastic lecturer. His popular lectures on mathematics and science designed for a wide audience sometimes attracted hundreds listeners. He had a lot of tricks to make his lectures very lively, attractive, unforgettable. For example, in 1940s he loved to give a talk "The Zhukovskii theorem on the supporting force of a plane's wing". The presentation, which was aimed at a very wide audience,  started as follows:

\noindent
\textit{Delone}: Do you know the famous mathematician, academician Paul Alexandrov?
\newline \textit{Someone in the audience}: Yes, we do.
\newline \textit{Delone}: Then, do you know that the academician Paul Alexandrov never flies in a plane?
\newline \textit{Audience} is puzzled.
\newline \textit{Delone}: Do you have any idea why academician Alexandrov is afraid of flying?
\newline \textit{Audience} keeps puzzling
\newline \textit{Delone}: Academician Alexandrov never flies on plane because he does not know  the Zhukovski theorem. But I do and I fly very much. Today I will make you acquainted with this wonderful theorem and you also will never fear  flying.

At this point the audience was  ``warmed-up'' and ready to be acquainted with a ``wonderful theorem''.
%% say what the theorem is, or give a reference

We mention one more example, Delone's \textit{affine cat} which has firmly entered Russian mathematical folk lore. During his lectures on analytic geometry Delone drew with a chalk on a blackboard a cat and several its affine images. The affine cats looked funny and gave additional representation to the geometry of affine transformations.

Mathematics students of Moscow university write jokes about their most popular professors on the theme: Who and in what does a soup cook? A version of the joke relevant to Delone contained a bit of irony: "Delone used to cook a soup in an integer point lattice. True, the soup used to spill, but the \textit{obviousness} remained."

In 1934 Delone initiated and organized the  first  mathematical olympiad in Russia. Participants in this and subsequent olympiads have left very warm and amazing reminiscences of Delone. Talking to pupils, Delone  evaluated the olympiads  highly and loved to repeat that the only difference between an olympiad problem and a real research problem is that an olympiad problem takes few hours while the research problem takes few thousand hours. Later for many years B.N.  was often the head of olympiad organizing committee.

Enthusiasm for mountain climbing, as has been already said, started in his childhood and continued throughout his life; it was, likely, his strongest passion after mathematics. His admiration for the beauty of mountain peaks and ridges was the stronger for its geometrical associations. For Delone, one of the most attractive components of mountaineering was the psychological fine-tuning needed for overcoming difficulties faced in mountains, for gaining a final victory over a summit and over oneself. Delone possessed these special psychological
features and, there can be no doubt,  these  alpinist qualities helped him  gain victory over serious mathematics problems.

Delone's mountain climbing was quite professional: he held a title, \textit{ master of Soviet mountaineering},
 organized several alpinist camps, and served  as a climbing coach. He was one of the best experts on the Western Caucasus. In 1937 he published a guide-book for alpinists "Summits of the Western Caucasus". This book contains detailed descriptions of paths on all the main mountains of this region. All of the numerous sketch-maps of routes of ascent on the tens of summits were drawn by Delone himself.
But the main feature of the guide-book is a very detailed panorama of the ridge chain of the whole Western Caucasus (over 200 kilometers long).  This panorama of the ridge chain was drawn by Delone himself sequentially from several outlying mountains of the region.

He also knew the Altai very well, especially the region of the Belukha. One of the paths up to the Belukha goes along an edge in which there is an intermediate mountain named the Delone Peak in the honor of one of its first subjugators.
However, his most beloved place in the Altai was not Delone Peak  but the lake Shavlo under the Skazka mountain  (``skazka'' is the Russian word for ``tale''). He visited this beautiful place many times, the last time in 1970, in his jubilee year. Delone was convinced that the beauty of this place exceeded any place in the Alps (he knew the Swiss Alps very well) and admired  it  all his life.

\medskip
\textbf{``An Evening of Life''}

\smallskip
In 1970 the German academy Leopoldina sent its congratulations on the occasion of the 80th anniversary of its member Boris Delone and wished him "a nice and quiet evening of life". The evening  continued for 10 more years and turned out to be very saturated, but not quiet at all.

Between his 80th and 90th years, Delone continued his very active life: research, popularization of science, care of family, hiking, etc..

   Delone was a very strong man from the birth and maintained his physical strength throughout his life by  weekly hiking, climbing, etc..  Each Saturday and Sunday he went on a 30-kilometer walking or skiing tour, usually a very remote path going through wild but beautiful places. These weekly tours took place in any weather, even very bad and ugly. Usually  a group of 5 or 6 people joined him in these tours and enjoyed the beauty of nature. Delone's intercourse with nature was sacramental. His close colleagues, who knew him well, recognized those instants when one should keep silence in order not to violate the voice of nature.

   On Summer holidays Delone went to real mountains: the Caucasus, the Altai, Tien Shan, the Pamirs. These tours were filled with varied, incredible events. Stories of his climbing or mountain hiking tours told by Delone himself were of extraordinary interest and could constitute the whole book.  Here is just one small typical episode which occurred to Delone in his 86th year.

July 6, 1975, Delone spends a cold night ( -25C) in a tent on a glacier under the famous peak of Khan-Tengri (7000 m, Tien-Shan mountain system, Central Asia) at a height about 4200 m. In the morning he comes down by helicopter to Frundze (now Bishkek, the capital of Kirgizstan) where heat exceeds +40C. After standing  few hours in line he succeeds in purchasing an air ticket to Moscow and arrived that night in a sub-Moscow airport. Taking the last local train he arrives at the very small station of Abramtsevo at 2 am and walks in deep night through a  dark dense forest to his "dacha" (country house). He loses the way and leaves his heavy rucksack in a secluded place. Only in the morning Delone succeeds in getting home safely.

It is not surprising that Delone was held in respect and liked by all his family. In his turn he remained a real head and support for his family and sought to resolve many of their problems. In particular, he cared very much for his beloved grandson Vadim Delone. The young Vadim was a talented poet and known dissident. On August 25, 1968, he took part in a celebrated political demonstration organized by 6 dissidents on the Red Square in Moscow against the Soviet Army's intervention in Czechoslovakia. For "one minute of freedom" with an unfolded slogan "For your and our freedom," the 20 year old Vadim was sentenced to 3 years imprisonment in a Siberian concentration camp. Delone did his best to ease his grandson's fate. He wrote letters to the KGB and asked many outstanding people to show support for Vadim, but nothing  helped. All he was able to do  was to visit the grandson in the camp as many times as he was allowed. One should add that Delone was not alone: Delone's former student Alexander Alexandrov\footnote{Academician A.D. Alexandrov (1912 - 1999) is an outstanding geometer of the XX century.} was also able  to get permission from the authorities and joined Delone in visiting Vadim.

There was one more recipient of his direct care. This was his wife Mariya Henrikhovna. She originated from a Denmark immigrant family. Since 1924 when they married, Mariya was B.N.'s closest and hopeful friend.
%% What is a hopeful friend?
They were born on the same day, March 15, though she was one year older.  That Delone had a long, successful, highly creative life was in great part due to her. In particular She made sure that he would not miss his walking, hiking, and climbing tours. One day  she told me  how  happy she was that B. N. had young colleagues and friends who were ready to join his weekly tours. When, at age about 85, she became seriously ill, Delone took care of her  to the end. This care was very touching. For two or three years, these duties seriously limited his possibilities for walking. Nevertheless, even at this time he had at least one serious walking tour (30-40 kilometers long) a week.

After his wife's death in February, 1976, Delone's "one hundred days" began. He through himself into a storm of activity: daily research, lecturing, and  on his two days off, 30 kilometers of hiking. Sometimes, on a working day that did not go well, B.N. made a suggestion: "Let us have some tea". This meant that he briefly stopped at  home to get his well-known campfire teapot. After that we, two or three, went to a rail-station (1 hour), took a train (over 1 hour), then walked in deep snow (half an hour) until we reached a beautiful glade surrounded by high splendid pine trees. We cooked tea made of snow, enjoyed a very nice sunny day for about an hour and started back, in order to be home in the metropolis by the late evening. On the train B.N. became an inexhaustible storyteller, another of his many talents. "Do you know that my cousin Liza Pilenko is holy?" he might suddenly ask. "What
 do you mean?" ask his puzzled companions (fellow-travelers.
 "Yes, yes, true, really holy", B. N. insists with a sly smile. After a theatrical pause he starts telling a  long, amazing, story about his cousin. After the revolution in 1917, a young poetess Liza Pilenko (married name Kuz'mina-Karavaeva) emigrated to Paris, where she changed professions several times, took the veil, and took the name Mother Mariya. In the early 1940's she entered the French Resistance, was imprisoned, and died in a concentration camp. Later the Catholic church officially canonized her.

One week in May turned out too stormy and ended with a hemorrhage of the brain. The "one hundred days" were over, "the Waterloo battle" had begun. Very limited mobility and speech, loss of  memory, age 86:  this promised nothing good. The doctors refused to give any consolatory prognosis.
Nevertheless, B.N. started fighting for his health, for a life of full value. Certainly, he was not alone in this battle. His relatives, students, colleagues did their best in help. And B.N. won his Waterloo battle: indeed, a few months later he returned to his normal life. Certainly, his schedule became less hectic: instead of 30 kilometer walking tours now he did just 10. Wild, impassable, remote forests gave  way to more "industrial" landscapes close to Moscow. However, B.N. had  a life of full value. Though this victory over such serious illness at such an age did not bring him additional fame,  it was likely the most valuable victory in his life.

Nevertheless, his life approached its natural end. I will never forget, how in March 1980, three or four months before he died, he took a book and read one remarkable statement by Poincar\'e. Overcome by emotion, his voice was trembling: "Life is just a short episode between two eternities of death and ... even in this episode conscious thought has lasted and will last only an instant.  A thought is just a lightning flash in the middle of the infinite night. But this flash is everything".

\end{document}